\documentclass{article}

\usepackage{arxiv}

\usepackage[utf8]{inputenc} % allow utf-8 input
\usepackage[T1]{fontenc}    % use 8-bit T1 fonts
\usepackage{hyperref}       % hyperlinks
\usepackage{url}            % simple URL typesetting
\usepackage{booktabs}       % professional-quality tables
\usepackage{amsfonts}       % blackboard math symbols
\usepackage{nicefrac}       % compact symbols for 1/2, etc.
\usepackage{microtype}      % microtypography
\usepackage{lipsum}		% Can be removed after putting your text content
\usepackage{graphicx}
\usepackage{natbib}
\usepackage{doi}

\usepackage[english]{babel}
\usepackage[version=4]{mhchem}
\usepackage{mathtools}
\usepackage{siunitx}
\usepackage{csvsimple}
\usepackage{verbatim}
\usepackage{tabularx} 

\newcommand{\B}[1]{\boldsymbol{#1}} % bold symbols
\newcommand{\R}{\mathbb{R}}
\newcommand{\x}{\times}
\newcommand{\del}{\partial}
\newcommand{\D}{\mathrm{d}}
\newcommand\norm[1]{\left\Vert#1\right\Vert}

\title{A Generalized Variable Projection Algorithm for Least Squares Problems in Atmospheric Remote Sensing}

\author{ \href{https://orcid.org/0009-0008-5497-1941}{\includegraphics[scale=0.06]{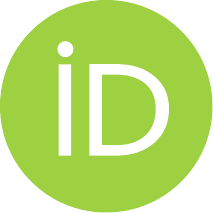}\hspace{1mm}Adelina Bärligea} \\
	Remote Sensing Technology Institute\\
	German Aerospace Center\\
	82234 Oberpfaffenhofen\\
	\And
	\href{https://orcid.org/0000-0001-9537-3050}{\includegraphics[scale=0.06]{orcid.pdf}\hspace{1mm}Philipp Hochstaffl} \\
	Remote Sensing Technology Institute\\
	German Aerospace Center\\
	82234 Oberpfaffenhofen\\
	\And
	\href{https://orcid.org/0000-0001-7196-6599}{\includegraphics[scale=0.06]{orcid.pdf}\hspace{1mm}Franz Schreier} \\
	Remote Sensing Technology Institute\\
	German Aerospace Center\\
	82234 Oberpfaffenhofen\\
}

\date{}

\renewcommand{\headeright}{Preprint}
\renewcommand{\undertitle}{Preprint to \url{https://doi.org/10.3390/math11132839                     }\thanks{submitted to MDPI mathematics, received: 17 May 2023, accepted: 22 June 2023, published: 24 June 2023.} }
\renewcommand{\shorttitle}{\url{https://doi.org/10.3390/math11132839                     } }

\begin{document}
\maketitle

\begin{abstract}
This paper presents a solution for efficiently and accurately solving separable least squares problems with multiple datasets. These problems involve determining linear parameters that are specific to each dataset while ensuring that the nonlinear parameters remain consistent across all datasets. A well-established approach for solving such problems is the variable projection algorithm introduced by Golub and LeVeque, which effectively reduces a separable problem to its nonlinear component. However, this algorithm assumes that the datasets have equal sizes and identical auxiliary model parameters. This article is motivated by a real-world remote sensing application where these assumptions do not apply. Consequently, we propose a generalized algorithm that extends the original theory to overcome these limitations. The new algorithm has been implemented and tested using both synthetic and real satellite data for atmospheric carbon dioxide retrievals. It has also been compared to conventional state-of-the-art solvers, and its advantages are thoroughly discussed. The experimental results demonstrate that the proposed algorithm significantly outperforms all other methods in terms of computation time, while maintaining comparable accuracy and stability. Hence, this novel method can have a positive impact on future applications in remote sensing and could be valuable for other scientific fitting problems with similar properties.
\end{abstract}

% keywords can be removed

\section{Introduction}
\label{sec:Intro}

One of the fundamental tasks in scientific computing is to find the unknown vector of parameters $\B{x}\in\R^{n+p}$, which, for given data $(y_{i},t_{i})$, $i=1,\ldots m$, $m\geq n+p$, and model function $\eta(\mathbf{x},t)$, solves the least squares problem
\begin{equation}
\min_{\B{x}}\sum_{i=1}^{m}(y_{i}-\eta(\B{x},t_{i}))^{2}=\min_{\B{x}}\bigl\Vert\B{y}-\B{\eta}(\B{x})\bigr\Vert_{2}^{2},
\label{eq:INnonlinearLeastSquaresProblem}
\end{equation}
where the vectors $\B{\eta}\in\R^{m}$ and $\B{y}\in\R^{m}$ have the entries\mbox{~$[\B{\eta}(\B{x})]_{i}=\eta(\B{x},t_{i})$} and \mbox{$[\B{y}]_{i}=y_{i}$,} respectively. If the model function is nonlinear, the minimization problem is solved iteratively by using step-length or trust-region methods, such as the Levenberg--Marquardt algorithm \citep{More.1978}. 

In separable least squares problems, the model function $\eta$ is a linear combination of nonlinear functions $\varphi_{j}$, $j=1,\ldots,n$, i.e.,
\begin{equation}
\eta(\B{\alpha},\B{\beta},t)=\sum_{j=1}^{n}\beta_{j}\varphi_{j}(\B{\alpha},t),
\label{eq:INseparablemodelfunction}
\end{equation}
and the vector of parameters $\B{x}$ is split into a vector of linear parameters $\B{\beta}\in\R^{n}$, and a vector of nonlinear parameters, $\B{\alpha}\in\R^{p}$. 

In 1973, \citet{Golub.1973} proposed a powerful method for solving such problems, called variable projection (VP). Specifically, the problem is reduced to a nonlinear least squares problem involving $\B{\alpha}$ only, and a linear least squares problem involving $\B{\beta}$ only. Thus, the dimension of the nonlinear problem to be solved is reduced from $n+p$ to $p$. Later on, \citet{G.H.Golub.1979} extended the method to the case of multiple datasets $\B{y}_{k}$ with $k=1,\ldots,s$. In this type of problem, a vector of linear parameters $\B{\beta}_{k}$ is associated with a dataset $\mathbf{y}_{k}$, while the vector of nonlinear parameters $\B{\alpha}$ corresponds to all datasets.

One example of this type of problem can be found in the retrieval of the atmospheric composition from passive remote sensing measurements (cf.\ \citep{Hochstaffl}). Here, an atmospheric radiative transfer model with molecular composition parameters (nonlinear) and reflectivity parameters (linear) is fitted to spectral radiance measurements to determine the underlying atmospheric state parameters. For long-lived molecules, such as \ce{CO2} or \ce{CH4}, which are homogeneously spread throughout the atmosphere, it is possible to fit several (nearby) observations for one concentration value while the surface reflectivity differs for each measurement, making this a separable problem with multiple right-hand sides.

Fitting, for example, $\num{4}\x\num{4}=\num{16}$ radiance spectra with three distinct linear reflectivity variables simultaneously for two nonlinear molecular concentration variables, the total number of unknowns for conventional solvers adds up to \num{50} ($=\num{16}\cdot\num{3}+\num{2}$). However, with a separable approach, such as variable projection, one can reduce this size by a factor of \num{25}. 

However, the variable projection algorithm established by \citet{G.H.Golub.1979} to solve such problems with multiple right-hand sides is based on the assumption that all datasets have the same lengths and corresponding nonlinear models. 

In this article, the theory by \citet{G.H.Golub.1979} is, therefore, elaborated and further enhanced for scenarios such as differently sized datasets or varying nonlinear model setups (see Section  \ref{sec:Theory}). The modifications were not only necessary in order to apply their method to the example described above, but this new algorithm could also be useful for other scientific fitting problems. The method was implemented in \textit{Python} (see Section \ref{sec:Algorithm}) and applied to the least squares problem arising in atmospheric trace gas retrieval (see Section \ref{sec:Experiments}). In Section \ref{sec:Conclusions}, experimental results are discussed and concluded. 

\section{Theoretical Background}
\label{sec:Theory}

Assuming a model that is a linear combination of nonlinear functions, such as \eqref{eq:INseparablemodelfunction}, one can define the arising separable least squares problem as follows: 
\begin{equation}
\min_{\B{\alpha},\B{\beta}}\sum_{i=1}^{m}\Bigl[y_{i}-\sum_{j=1}^{n}\beta_{j}\varphi_{j}(\B{\alpha},t_{i})\Bigr]^{2}=\min_{\B{\alpha},\B{\beta}}\bigl\Vert\B{y}-\B{\Phi}(\B{\alpha})\B{\beta}\bigr\Vert_{2}^{2},\label{eq:INseparableFunctional}
\end{equation}
where the matrix $\B{\Phi}(\B{\alpha})\in\R^{m\x n}$
has the entries $[\B{\Phi}(\B{\alpha})]_{ij}=\varphi_{j}(\B{\alpha},t_{i})$.

\subsection{The Variable Projection Method}
\label{subsec:VP}

The minimization problem in (\ref{eq:INseparablemodelfunction}), written as 
\begin{equation}
\min_{\B{\alpha}}\Bigl(\min_{\B{\beta}}\bigl\Vert\B{y}-\B{\Phi}(\B{\alpha})\B{\beta}\bigr\Vert_{2}^{2}\Bigr),
\label{eq:THtwoleastsquares}
\end{equation}
is reduced to two least squares problems. \citet{Golub.1973} pointed out that for any fixed $\B{\alpha}$, one considers the merely linear least squares problem 
\begin{equation}
\min_{\B{\beta}}\bigl\Vert\B{y}-\B{\Phi}(\B{\alpha})\B{\beta}\bigr\Vert_{2}^{2}.
\label{eq:THlinearLS}
\end{equation}
which is solved by
\begin{equation}
\B{\beta}(\B{\alpha})=\B{\Phi}^{\dagger}(\B{\alpha})\B{y},
\label{eq:INlinearleastsquaressolution}
\end{equation}
where $\B{\Phi}^{\dagger}=(\B{\Phi}^{\top}\B{\Phi})^{-1}\B{\Phi}^{\top}$ is the generalized inverse of the matrix $\B{\Phi}$. Second, for the~functional
\begin{equation}
\mathcal{F}(\B{\alpha})=\min_{\B{\beta}}\bigl\Vert\B{y}-\B{\Phi}(\B{\alpha})\B{\beta}\bigr\Vert_{2}^{2}=\bigl\Vert\B{y}-\B{\Phi}(\B{\alpha})\B{\Phi}^{\dagger}(\B{\alpha})\B{y}\bigr\Vert_{2}^{2},
\label{eq:INvariableprojectionfunctional}
\end{equation}
which is obtained by substituting the expression of $\B{\beta}(\B{\alpha})$ given by Equation (\ref{eq:INlinearleastsquaressolution}) into the residual function, one considers the nonlinear least squares problem
\begin{equation}
\min_{\B{\alpha}}\mathcal{F}(\B{\alpha})=\min_{\B{\alpha}}\bigl\Vert\B{y}-\B{\Phi}(\B{\alpha})\B{\Phi}^{\dagger}(\B{\alpha})\B{y}\bigr\Vert_{2}^{2}=\min_{\B{\alpha}}\bigl\Vert\B{P}_{\B{\Phi}(\B{\alpha})}^{\perp}\B{y}\bigr\Vert_{2}^{2},
\label{eq:INvariableprojectionfunctional2}
\end{equation}
where $\B{P}_{\B{\Phi}(\B{\alpha})}^{\perp}=\B{I}_{m\x m}-\B{P}_{\B{\Phi}(\B{\alpha})}\in\R^{m\x m},$ and $\B{P}_{\B{\Phi}}(\B{\alpha})=\B{\Phi}(\B{\alpha})\B{\Phi}^{\dagger}(\B{\alpha})\in\R^{m\x m}$ is the orthogonal projection operator onto the column space of the matrix $\B{\Phi}(\B{\alpha})$. Thanks to this formulation, the method is called variable projection (VP). 

Thus, the method of solution involves, in the first step, the computation of the minimizer $\widehat{\B{\alpha}}=\arg\min_{\B{\alpha}}\mathcal{F}(\B{\alpha})$, and in the second step, the computation of the vector of linear parameters as $\widehat{\B{\beta}}=\B{\Phi}^{\dagger}(\widehat{\B{\alpha}})\B{y}$. \citet{Golub.1973} showed that this solution method yields the correct result, under the assumption that $\B{\Phi}(\B{\alpha})$ has a locally constant rank in the neighborhood of $\widehat{\B{\alpha}}$. 
In order to minimize a nonlinear problem, such as \eqref{eq:INvariableprojectionfunctional}, the derivative of the objective function, with respect to the unknown variables, is needed. If the rank of matrix $\B{\Phi}(\B{\alpha})$ was not constant across the points throughout the iteration where its derivative is calculated, the pseudo-inverse $\B{\Phi^{\dagger}}(\B{\alpha})$ would not be a continuous function and, therefore, not differentiable.

This separation technique has powerful advantages over conventional algorithms, which all follow from the fact that the nonlinear functional $\mathcal{F}(\B{\alpha})$ only depends on $\B{\alpha}$ (see also the review by \citet{Golub.2003}). Since $\B{\alpha}$ is a vector of length $p$ and $\B{\beta}$ of length $n$, the VP method effectively turns the original nonlinear problem of $n+p$ variables into one of only $p$. This reduction of the parameter space of the problem results in a smaller size of the problem's Jacobian and, therefore, requires less time for its computation. Moreover, \citet{OLeary.2013} pointed out that the reduced parameter space may also lead to a reduced number of local minimizers, making it more likely to find the global minimum instead of a local one. 
To summarize, the VP solver is a lot more efficient and also converges better than conventional methods with no separation \citep{Ruhe.1980b, Golub.2003}. Moreover, the need for a smaller initial guess vector can generally lead to a better-conditioned and more stable problem. 

For solving the nonlinear least squares problem (\ref{eq:INvariableprojectionfunctional2}), one needs to calculate the partial derivatives $\del\B{P}_{\B{\Phi}}^{\perp}/\del\alpha_{l}$. In this regard, \citet{Golub.1973} provide the {computational~formula} 
\begin{equation}
\frac{\del\B{P}_{\B{\Phi}}^{\perp}}{\del\alpha_{l}}=-\Bigl[\B{P}_{\B{\Phi}}^{\perp}\frac{\del\B{\Phi}}{\del\alpha_{l}}\B{\Phi}^{\dagger}+\Bigl(\B{P}_{\B{\Phi}}^{\perp}\frac{\del\B{\Phi}}{\del\alpha_{l}}\B{\Phi}^{\dagger}\Bigr)^{\top}\Bigr],
\label{eq:INjacobianofVPfunctional}
\end{equation}
which is proved in the Appendix. 

\citet{Kaufman.1975} proposed a modification of the original method by making use of a QR decomposition of the matrix $\B{\Phi}$:
\begin{equation}
\B{\Phi}=\B{Q}\B{R}=(\begin{array}{cc} 
\B{Q}_{1} & \B{Q}_{2})\left(\begin{array}{c}
\B{R}_{1}\\
\B{0}_{(m-n)\x n}
\end{array}\right)=\B{Q}_{1}\B{R}_{1},\end{array}
\label{eq:INqrdecomposition}
\end{equation}
where $\B{Q}\in\R^{m\x m}$, $\B{Q}_{1}\in\R^{m\x n}$, and $\B{Q}_{2}\in\R^{m\x(m-n)}$ are matrices with orthonormal columns (i.e., $\B{Q}^{\top}\B{Q}=\B{I}_{m\x m}$, $\B{Q}_{1}^{\top}\B{Q}_{1}=\B{I}_{n\x n}$, $\B{Q}_{2}^{\top}\B{Q}_{2}=\B{I}_{(m-n)\x(m-n)}$) and $\B{R}_{1}\in\R^{n\x n}$ is an upper triangular, non-singular matrix. In this case, the generalized inverse $\B{\Phi}^{\dagger}$, satisfying the relation $\B{\Phi}^{\dagger}\B{\Phi}=\B{I}_{n\x n}$, computes as 
\begin{equation}
\B{\Phi}^{\dagger}=(\begin{array}{cc}
\B{R}_{1}^{-1} & \B{0}_{n\x(m-n)})\left(\begin{array}{c}
\B{Q}_{1}^{\top}\\
\B{Q}_{2}^{T}
\end{array}\right)=\B{R}_{1}^{-1}\B{Q}_{1}^{\top}.\end{array}
\label{eq:THpseudoinverseQR}
\end{equation}

From the invariance of the $2$-norm under orthogonal transformations and  $\B{Q}^{\top}\B{\Phi}=\B{Q}^{\top}\B{Q}\B{R}=\B{R}$), it follows that
\begin{align}
\mathcal{F}(\B{\alpha}) & =\min_{\B{\beta}}\bigl\Vert\B{y}-\B{\Phi}(\B{\alpha})\B{\beta}\bigr\Vert_{2}^{2}\nonumber \\
 & =\min_{\B{\beta}}\bigl\Vert\B{Q}^{\top}(\B{\alpha})(\B{y}-\B{\Phi}(\B{\alpha})\B{\beta})\bigr\Vert_{2}^{2}\nonumber \\
 & =\min_{\B{\beta}}\Bigl\Vert\left(\begin{array}{c}
\B{Q}_{1}^{\top}(\B{\alpha})\B{y}\\
\B{Q}_{2}^{\top}(\B{\alpha})\B{y}
\end{array}\right)-\left(\begin{array}{c}
\B{R}_{1}(\B{\alpha})\\
\B{0}_{(m-n)\x n}
\end{array}\right)\B{\beta}\Bigr\Vert_{2}^{2}\nonumber \\
 & =\min_{\B{\beta}}\bigl\Vert\B{Q}_{1}^{\top}(\B{\alpha})\B{y}-\B{R}_{1}(\B{\alpha})\B{\beta}\bigr\Vert_{2}^{2}+\bigl\Vert\B{Q}_{2}^{\top}(\B{\alpha})\B{y}\bigr\Vert_{2}^{2}.
\label{eq:INqrPrewritten}
\end{align}

Since the optimal $\B{\beta}$ for any given $\B{\alpha}$ is
\begin{equation}
    \B{\beta}(\B{\alpha})=\B{\Phi}^{\dagger}(\B{\alpha})\B{y}=\B{R}_{1}^{-1}(\B{\alpha})\B{Q}_{1}^{\top}(\B{\alpha})\B{y},
\end{equation}
we find $\mathcal{F}(\B{\alpha})=\bigl\Vert\B{Q}_{2}^{T}(\B{\alpha})\B{y}\bigr\Vert_{2}^{2}$, showing that the nonlinear least squares problem reduces to
\begin{equation}
\min_{\B{\alpha}}\mathcal{F}(\B{\alpha})=\min_{\B{\alpha}}\bigl\Vert\B{Q}_{2}^{\top}(\B{\alpha})\B{y}\bigr\Vert_{2}^{2}.
\label{eq:INqrrewrittenfrom}
\end{equation}

Moreover, for derivative calculations, \citet{Kaufman.1975} proposed the simplified formula 
\begin{equation}
\frac{\del\B{Q}_{2}^{T}}{\del\alpha_{l}}=-\B{Q}_{2}^{\top}\frac{\del\B{\Phi}}{\del\alpha_{l}}\B{\Phi}^{\dagger},
\label{eq:THderivativeQ2}
\end{equation}
which is justified in the Appendix. This was shown to save function and gradient evaluation costs and, therefore, reduce the computing time per iteration, which is why this simplified version of \eqref{eq:INjacobianofVPfunctional} became well established.

In conclusion, the minimization problems (\ref{eq:INvariableprojectionfunctional2}) and (\ref{eq:INqrrewrittenfrom}) are equivalent, but the size of the matrix $\B{Q}_{2}^{\top}$ is smaller than that of the matrix $\B{P}_{\B{\Phi}}^{\perp}$. Consequently, an algorithm for solving  Equation (\ref{eq:INqrrewrittenfrom}) should be more efficient. 

\subsection{Multiple Right-Hand Sides (MRHS)}
\label{subsec:Golub_leVeque Approach}

As \citet{G.H.Golub.1979} pointed out, there can be optimization problems where multiple sets of data $\B{y_1}, \B{y_2}, \ldots, \B{y_s}$ are to be fit to a model function, such as \eqref{eq:INseparablemodelfunction}. If the model parameters are to vary for each set $\B{y_k}$, this will just result in $s$ distinct separable problems, such as \eqref{eq:INseparableFunctional}. There are, however, cases where only the linear variables are specific to each dataset, while the nonlinear variables have to hold for all available data simultaneously. By exploiting separability, the sizes of such minimization problems can be reduced from $sn+p$ to just $p$ unknown variables, which is even greater than with only one dataset.

Separable problems with $s$ right-hand side(s) (RHS) can be posed as the minimization of
\begin{equation}
    \min_{\B{\alpha}, \B{B}}\ \norm{\B{Y}-\B{\Phi}(\B{\alpha})\B{B}}_\text{F}^2,
\label{eq:THmatrixMINforMRHS}
\end{equation}
where the nonlinear parameter vector is $\B{\alpha}\in \R^p$ and the matrix is $\B{B} \in \R ^{n\x s}$, containing the linear parameter vectors $\B{\beta_k}$ for each RHS $\B{y_k}$ with $k=1,\ldots,s$ as its columns, using the Frobenius norm $\norm{\cdot}_\text{F}$. Here, the data matrix $\B{Y}=(\B{y_1} \hdots \B{y_s}) \in\R ^{m\x s}$ is fit to a single model of the form $\B{\Phi}(\B{\alpha})\B{B}$, where $\B{\Phi}(\B{\alpha})\in\R^{m\x n}$ is as before (cf.\ Equation\ \eqref{eq:INseparableFunctional}). 

\subsubsection{Naive Approach}
\label{ssec:NAIVE}
The first intuitive approach, which was also mentioned by \citet{G.H.Golub.1979}, is to reformulate \eqref{eq:THmatrixMINforMRHS} using the matrix 
\begin{equation}
\B{G}(\B{\alpha}) = 
    \begin{pmatrix}
\B{\Phi}(\B{\alpha}) &  & \hdots & \B{0} \\
 & \B{\Phi}(\B{\alpha}) & & \vdots \\
 \vdots & & \ddots & \\
 \B{0} & \hdots & & \B{\Phi}(\B{\alpha})
     \end{pmatrix}
\in \R^{ms\x ns},
\label{eq:THzeromatrix}
\end{equation}
and the vectors $ \B{\Tilde{y}} = 
\begin{pmatrix} \B{y_1} \\ \vdots  \\ \B{y_s} \end{pmatrix} \in \R^{ms}$ and 
$\B{\Tilde{\beta}} = 
\begin{pmatrix} \B{\beta_1} \\ \vdots \\ \B{\beta_s} \end{pmatrix} \in \R^{ns}$, 
such that a problem of the original vectorial form
\begin{equation}
\min_{\B{\alpha}, \B{\Tilde{\beta}}}\ \norm{\B{\Tilde{y}}-\B{G}(\B{\alpha})\B{\Tilde{\beta}}}^2,
\label{eq:THnaiveproblem}
\end{equation}
arises. This formulation allows for solving the separable problem with multiple RHS, by means of the variable projection algorithm already discussed. However, even for a moderate number of datasets, the matrix $\B{G}(\B{\alpha})$ becomes overly large. Moreover, the sparse structure and the fact that all diagonal blocks in \eqref{eq:THzeromatrix} are the same can be better utilized in the earlier formulation  \eqref{eq:THmatrixMINforMRHS}. 

\subsubsection{Golub--LeVeque Approach}
\label{ssec:GOLUB}
Therefore, \citet{G.H.Golub.1979} suggested a different approach: Starting from \eqref{eq:THmatrixMINforMRHS}, one can, in the same manner as Equation \eqref{eq:INvariableprojectionfunctional}, exploit the problem's separability and reduce it to a purely nonlinear minimization problem of the form 
\begin{equation}
    \min_{\B{\alpha}}\ \norm{\B{P}^{\perp}_{\B{\Phi}}(\B{\alpha})\B{Y}}_\text{F}^2,
\label{eq:THreducedProblemMATRIX}
\end{equation}
with the orthogonal projector $\B{P}^{\perp}_{\B{\Phi}}(\B{\alpha}) = \B{I}-\B{\Phi}(\B{\alpha})\B{\Phi^{\dagger}}(\B{\alpha})$ as before. Acknowledging that the Frobenius norm in \eqref{eq:THreducedProblemMATRIX} is equivalent to the \num{2}-norm of a vector function $\B{z}(\B{\alpha})$, set up as 
\begin{equation}
    \B{z}(\B{\alpha}) = \begin{pmatrix}
 \B{P}^{\perp}_{\B{\Phi}}(\B{\alpha})\: \B{y_1}\\
  \vdots \\ \B{P}^{\perp}_{\B{\Phi}}(\B{\alpha})\: \B{y_s} \end{pmatrix} \in \R^{ms},
\label{eq:THz(a)}
\end{equation}
this can be minimized with any of the established methods for nonlinear least squares problems, such as the Levenberg--Marquardt algorithm, which is used in \citep{G.H.Golub.1979}. The Jacobian matrix of $\B{z}(\B{\alpha})$ can be calculated analogously by defining its $l{\text{th}}$ column as 
\begin{equation}
    \frac{\del \B{z}(\B{\alpha})}{\del \alpha_l} = \begin{pmatrix}
 \frac{\del \B{P}^{\perp}_{\B{\Phi}}(\B{\alpha})}{\del \alpha_l}\: \B{y_1}\\
  \vdots \\ \frac{\del \B{P}^{\perp}_{\B{\Phi}}(\B{\alpha})}{\del \alpha_l}\: \B{y_s} \end{pmatrix} \in \R^{ms},
\label{eq:THdz(a)/dal}
\end{equation}
where $\frac{\del \B{P}^{\perp}_{\B{\Phi}}(\B{\alpha})}{\del \alpha_l}$ is exactly as in \eqref{eq:INjacobianofVPfunctional}. 

\subsubsection{Kaufman Approach}
\label{ssec:KAUFMAN}
Since the Frobenius norm is invariant under the orthogonal transformation, the above method can also be written in terms of the QR decomposition of $\B{\Phi}(\B{\alpha})$ outlined in \mbox{Section  \ref{subsec:VP}.} This approach based on \citet{Kaufman.1975} can, therefore, be established by replacing the orthogonal projector $\B{P}^{\perp}_{\B{\Phi}}(\B{\alpha})\in \R^{m\x m}$ in all of the above equations of the Golub--LeVeque method, with the smaller orthogonal matrix $\B{Q_2}^{\top}\in \R^{(m-n)\x m}$ derived in \eqref{eq:INqrdecomposition}. Both versions of the approach are included in the \textit{Python} implementation outlined in Section \ref{sec:Algorithm} and are, therefore, the subjects of the numerical experiments performed in Section \ref{sec:Experiments}. 

\subsection{Extensions to the Golub--LeVeque Approach}
\label{subsec:Extentions to GL}

Working with real measurements, which can be subject to errors or missing data points, it is not always possible to use datasets that all have the exact same number of data points. In this case, each dataset $\B{y_k}$ has a specific size $m_k$, such that the matrix depiction used \mbox{in \eqref{eq:THmatrixMINforMRHS}} does not hold any longer.

Moreover, the nonlinear model functions stored in the matrix $\B{\Phi}(\B{\alpha})$ often depend on further auxiliary parameters, which may vary for each individual measurement (e.g.,\ observation angle). Thus, there can be different $\B{\Phi_k}(\B{\alpha})$ of size $m_k\x n$.  

None of these aspects were explicitly mentioned by \citet{G.H.Golub.1979}
or \citet{Kaufman.1992}, who later simplified the Golub--LeVeque method further,  with the assumption of equal lengths of the $\B{y_k}$. It was only mentioned by \citet{LindaKaufman.2010} that the TIMP package by \citet{Mullen2007} allows for differently sized data vectors, but not for differently constructed $\B{\Phi}$s. 

For the implementation introduced in this paper, both of these changes are taken into account. It can be seen that the special structure of vector \eqref{eq:THz(a)} and matrix \eqref{eq:THdz(a)/dal} can be exploited to rewrite 
\begin{equation}
    \B{z}(\B{\alpha}) = \begin{pmatrix}
 \B{P}^{\perp}_{\B{\Phi_1}}(\B{\alpha})\: \B{y_1}\\
  \vdots \\ \B{P}^{\perp}_{\B{\Phi_s}}(\B{\alpha})\: \B{y_s} \end{pmatrix} \in \R^{m_1+\ldots+m_s}
\label{eq:THz(a)NEW}
\end{equation}
as a stack of (differently sized) vectors and, likewise, the $l{\text{th}}$ column of its Jacobian
\begin{equation}
    \frac{\del \B{z}(\B{\alpha})}{\del \alpha_l} = \begin{pmatrix}
 \frac{\del \B{P}^{\perp}_{\B{\Phi_1}}(\B{\alpha})}{\del \alpha_l}\: \B{y_1}\\
  \vdots \\ \frac{\del \B{P}^{\perp}_{\B{\Phi_s}}(\B{\alpha})}{\del \alpha_l}\: \B{y_s} \end{pmatrix} \in \R^{(m_1+\ldots+m_s)s}\ \mathrm{.}
\label{eq:THdz(a)/dalnew}
\end{equation}
 
This necessary modification naturally increases the computational expense compared to the original Golub--LeVeque method, as matrix $\B{\Phi}$ and its derivative have to be calculated $s$ times instead of once. This also reduces the gain in efficiency that one would have had from exchanging $\B{P}^{\perp}_{\B{\Phi}}$ with $\B{Q_2}^{\top}$, as the QR decomposition now needs to be calculated for every single matrix $\B{\Phi}_k$.
However, this modified VP method for multiple right-hand sides is still significantly more efficient than using a standard nonlinear optimizer for the same problem, as will be shown in Section \ref{sec:Experiments}.

\section{Implementation}
\label{sec:Algorithm}
The first \textit{FORTRAN} implementation of a variable projection algorithm that allows for multiple right-hand sides was VARP2 \citep{NetlibRepository.09.05.1985.varpro2}, which was developed by Randy LeVeque. It is a modification of the subroutine VARPRO \citep{NetlibRepository.09.05.1985.varpro} for classical separable problems with only a single RHS. Both work with user-provided Jacobians. 
Another efficient implementation of a solver for multiple datasets is the TIMP package by \citet{Mullen2007}, written in the statistical computing language \textit{R}, where the Jacobian is calculated by finite \mbox{difference approximations. }

This work was partly motivated by the fact that many scientific computing problems are nowadays implemented in \textit{Python}. However, the only modern implementation of a VP algorithm is a \textit{MATLAB} code by \citet{OLeary.2013}, which does not allow for multiple right-hand sides, as discussed in Section  \ref{subsec:Extentions to GL}. Here, an implementation is introduced where the standard VP algorithm by \citet{OLeary.2013} is enhanced for multiple datasets, as described in Section \ref{sec:Theory}. 

The \textit{MATLAB} code by \citet{OLeary.2013} was specifically designed to be brief and easily understandable, so it would be well-suited for translation. \citet{OLeary.2013} argued that many of the established implementations, such as VARPRO \citep{NetlibRepository.09.05.1985.varpro}, and PORT \mbox{library \citep{Fox78}} subroutines, such as NSF/NSG \citep{Kaufman.1975}, are somewhat outdated today and lack readability; thus, they proposed an efficient present-day implementation written in an interpreted language, with the advantage that it can easily be enhanced or modified. One reason for their code's brevity is that they used built-in \textit{MATLAB} functions, such as \textit{lsqnonlin.m}, to solve the nonlinear least squares problem \eqref{eq:INvariableprojectionfunctional}, and \textit{svd.m} to solve the remaining linear problems via singular value decomposition, instead of writing their own, making the algorithm modular. Another advantage of their code is the variety of statistical diagnostics it offers, some of which were also modified for multiple datasets and used for the assessment in Section \ref{sec:Experiments}. 

\subsection{Algorithm for a Single Right-Hand Side}
\label{subsec:AlgorithmOutline}
A short outline of a standard VP algorithm (such as the one by \citet{OLeary.2013}) can be formulated as this:
\begin{enumerate}
    \item The user supplies the measurement vector $\B{y}$, the number of linear parameters $n$, additional independent variables for the model, and an initial guess for the nonlinear variables $\B{\alpha_0}$. In addition, a subroutine/function which will be called \textit{ADA} (note that a detailed description of how to set up and implement this subroutine for a variable projection algorithm, such as the one presented here, can be found in \citep{OLeary.2013}, Sections 2.3 and 2.4) needs to be provided, which, for the given input $n$ and $\B{\alpha}$, calculates the nonlinear model matrix $\B{\Phi}$ and its partial derivatives $\del \B{\Phi}/\del \alpha_l$. 
    \item Calculate the pseudo-inverse $\B{\Phi}^{\dagger}$ to generate the variable projection functional \mbox{from \eqref{eq:INvariableprojectionfunctional}}, which is dependent on a given $\B{\alpha}$.
    \item Use the partial derivatives $\del \B{\Phi}/\del \alpha_l$ to generate the Jacobian matrix $\frac{\D\B{P}^{\perp}_{\B{\Phi}}}{\D\B{\alpha}}\:\B{y}$ for a given $\B{\alpha}$, as in \eqref{eq:INjacobianofVPfunctional}.
    \item Minimize the variable projection functional \eqref{eq:INvariableprojectionfunctional2} by using the results from steps 2 and 3 and insert them into an already existing nonlinear least squares solver to obtain the final result $\hat{\B{\alpha}}$. 
    \item Take $\hat{\B{\alpha}}$ and calculate a final $\B{\Phi}^{\dagger}(\hat{\B{\alpha}})$ to solve the remaining linear parameters $\hat{\B{\beta}}$, as \mbox{in \eqref{eq:INlinearleastsquaressolution}.}
\end{enumerate}

In the original implementation by \citet{OLeary.2013}, the linear solution, including the computation of $\B{\Phi}^{\dagger}(\B{\alpha})$, was calculated via a singular value decomposition, which is a very robust direct method for solving linear least squares problems. Moreover, instead of using the simplification by \citet{Kaufman.1975} of dropping the second term in the Jacobian \eqref{eq:INjacobianofVPfunctional}, \citet{OLeary.2013} argued that, in modern computers, the balance between the computing time for extra iterations vs. the second term has tipped. In 1975, when \citet{Kaufman.1975} proposed  simplification, it used to be more efficient to put up with more functional evaluations by saving matrix computations; this could have changed today \citep{OLeary.2013}. To verify this argument, a simplification by \citet{Kaufman.1975} is included as a subject for testing in Section \ref{subsec:testsreal}. 

To solve the nonlinear problem in step 4, the \textit{\scalebox{.95}[1.0]{least\_squares}} function \citep{SciPyv1.8.0Manual.06.02.2022} from the \textit{scipy.optimize} package \citep{SciPy2020} is used in the implementation. Its default solving method is the so-called trust region reflective algorithm \citep{Branch.1999}, which has been shown to work efficiently within the VP algorithm for the example discussed in Section \ref{sec:Experiments} \citep{Barligea2022}. 

\subsection{Modifications for Multiple Right-Hand Sides}
\label{subsec:ModificationsforMRHS}
In the following, the modifications required for multiple right-hand sides $\B{y_1},\ldots,\B{y_s}$ with possibly varying lengths $m_1,\ldots,m_s$, are going to be discussed.
For the \textbf{naive approach}, the only necessary change of the above procedure arises in the user-supplied function \textit{ADA} from step 1; instead of one matrix $\B{\Phi}(\B{\alpha})$ and its derivative, it has to return the sparse matrix $\B{G}(\B{\alpha})$ as in \eqref{eq:THzeromatrix}, with the possibly different $\B{\Phi}_k(\B{\alpha})$ on its diagonal, and the corresponding derivatives are also sparse matrices with the Jacobians $\del \B{\Phi}_k/\del \B{\alpha}$ on the diagonal. After setting up $\B{\Tilde{y}}$ and $\B{\Tilde{\beta}}$, as in Equation \eqref{eq:THnaiveproblem}, the above steps can be performed unchanged to obtain the results $\hat{\B{\alpha}}$ and $\Tilde{\B{\beta}}=(\hat{\B{\beta_1}},\ldots,\hat{\B{\beta_s}})^{\top}$.

This is different from the \textbf{Golub--LeVeque approach} described in Section  \ref{subsec:Golub_leVeque Approach}. Here, the user-defined method \textit{ADA} remains the same as in step 1 of Section \ref{subsec:AlgorithmOutline}; however, instead of calculating \eqref{eq:INvariableprojectionfunctional2} in step 2, the new vector function $\B{z}(\B{\alpha})$ has to be set up, such as \eqref{eq:THz(a)NEW}. The same holds for the new derivative matrix $\D\B{z}(\B{\alpha})/\D\B{\alpha}$ from \eqref{eq:THdz(a)/dalnew} in step 3. After this, in step 4, $\norm{\B{z}(\B{\alpha})}^2$ is minimized for $\hat{\B{\alpha}}$; consequently, in step 5, the final linear parameter matrix $\hat{\B{B}}$ can be derived by calculating its $k{\text{th}}$ column as follows: 
\begin{equation}
    \B{\beta_k} = \B{\Phi_k}^{\dagger}(\hat{\B{\alpha}})\:\B{y_k}\ \mathrm{.}
\label{eq:ALfinalbetas}
\end{equation}

Finally, for the \textbf{Kaufman approach} \citep{Kaufman.1975}, it is important to note that the pseudo-inverse computed in steps 2 and 5 is now  calculated via the QR decomposition \eqref{eq:INqrdecomposition}, as in \eqref{eq:THpseudoinverseQR}, while the remaining matrix $\B{Q_2}$ of the decomposition can be taken for setting up $\B{z}(\B{\alpha})$ and $\del\B{z}(\B{\alpha})/\del\B{\alpha}$, just as before, by replacing $\B{P}^{\perp}_{\B{\Phi}}$ in \eqref{eq:THz(a)NEW} and \eqref{eq:THdz(a)/dalnew} with $\B{Q_2}^{\top}$. The simplified Jacobian \eqref{eq:THderivativeQ2} is also adopted in the implementation of this approach. 

\section{Numerical Experiments}
\label{sec:Experiments}
The goal of this section is to outline an example of a separable least squares problem with multiple right-hand sides and use it to evaluate the performance of the suggested VP algorithms. The tests were also set up to establish a comparison to conventional nonlinear least squares (NLS) methods, which ignore separability. Those were implemented using the \textit{least\_squares} function from the \textit{scipy.optimize} package \citep{SciPyv1.8.0Manual.06.02.2022}.

Without separation, the minimization problem of multiple right-hand sides can be composed as follows:
\begin{equation}
    \min_{\B{x}}\ \norm{\begin{pmatrix} \B{y_1} \\ \vdots \\ \B{y_s} \end{pmatrix} - \begin{pmatrix} \B{\eta}(\B{\beta_1},\B{\alpha}) \\ \vdots \\ \B{\eta}(\B{\beta_s},\B{\alpha}) \end{pmatrix} }^2 \qquad\text{with}\ \B{x}=(\B{\alpha},\B{\beta_1}, \ldots, \B{\beta_s})^{\top}\ \mathrm{.}
\label{eq:NEreferencemethod}
\end{equation}
which is easily extendable to the case of varying sizes of $\B{y_1},\ldots,\B{y_s}$, and correspondingly different $\B{\eta_1}, \ldots, \B{\eta_s}$. 

\subsection{Application: Trace Gas Retrieval}
A real-world example for the described problem set can arise in the area of remote sensing, more specifically, in the retrieval of atmospheric trace gas concentrations from spectral radiance measurements. The concentration of carbon dioxide (\ce{CO2}) or methane (\ce{CH4})---both important greenhouse gases---can be inferred from spectra observed in the short-wave infrared (SWIR). Such measurements are often spaceborne in order to achieve global coverage of the atmospheric composition. 

For this paper, observations from the OCO-2 (Orbiting Carbon Observatory-2) satellite by NASA \citep{CRISP2004700, Eldering2017, OCO2} were used, which was designed to monitor \ce{CO2} by measuring its absorption bands in the SWIR. In this spectral region, a radiance measurement can be modeled by the radiative transfer model (based on the Beer--Lambert law for molecular absorption neglecting scattering) \citep{GimenoGarcia.2011}
\begin{equation}
\hat{I}(\nu,\B{r},\B{\alpha}) =\ \left(\sum_{j=0}^{n-1} r_j\:\nu^j\right) \cdot \mu_{\odot} \cdot I_{0}(\nu) \cdot \exp\left(-\sum_{l=1}^p \alpha_l\: \tau_l^{\text{prior}}(\nu)\right )\ \otimes \ S(\nu),
\label{eq:NEradianceMODEL}
\end{equation}
dependent on the wavenumber $\nu$. The term wavenumber, which represents the inverse of the often used wavelength $\lambda$, has the common unit for the SWIR region of $\si{\cm\tothe{-1}}$, corresponding to $\num{e4}/\si{\micro m}$. The two sets of fitting parameters are $\B{r}\in\R^n$ for linear parameters $r_j$ and $\B{\alpha}\in\R^p$ for nonlinear parameters $\alpha_l$. 

The first term is a polynomial that approximates the wavenumber-dependent surface reflectivity of the Earth at the measurement location. Factor $\mu_{\odot}$ corresponds to $\cos(\theta_{\odot})$ (with $\theta_{\odot}$: solar zenith angle) and accounts for the geometry of the measurement, $I_{0}(\nu)$ is the incoming solar radiation (at the top of the atmosphere), and 
\begin{equation}
    \tau_l^{\text{prior}}(\nu) =\int_{s}n_l(s)\ \sigma_l(\nu,p(s),T(s))\ \D s
\label{eq:NEtotalopticaldepth1}
\end{equation}
is the total optical depth of the $l{\text{th}}$ molecule, which is the path integral over the number density $n_l$, and its pressure- and temperature-dependent cross-section $\sigma_l$. In trace gas retrieval, $\tau(\nu)$ is the most important measure, as it is directly related to a molecule's concentration in the atmosphere on a given path $s$. Unfortunately, SWIR observations do not provide enough information to retrieve the concentration profile of a molecule. Therefore, the calculation of \eqref{eq:NEtotalopticaldepth1} is only possible under prior assumptions of the atmospheric state (i.e.,\ the temperature and pressure profiles $p(s)$, $T(s)$, and the molecular number density, i.e., $n(s)$). Hence, a simple scaling factor $\alpha_l$ is fitted as
\begin{equation}
    \tau_l(\nu) = \alpha_l\cdot\tau_l^{\text{prior}}(\nu)
\label{eq:NEscalingTods}
\end{equation}
in the forward model \eqref{eq:NEradianceMODEL} to retrieve the ``real'' optical depth at the time and place of measurement. 
Lastly, $S(\nu)$ is the spectral response function of the sensor, which has to be convolved with the monochromatic radiance in order to mimic a real measurement.

For trace gas retrieval, one has to consider all $p$ molecules that have non-negligible absorbance in the measured spectral region. In the case of OCO-2 observations, the only relevant molecule apart from \ce{CO2} is \ce{H2O} (water vapor), meaning that there are two nonlinear fitting parameters. For the linear parameters, it is common to use approximately three reflectivity coefficients (depending on the size of the spectral interval). This means that for each spectrum, the necessary fitting parameters are
\begin{equation}
    \B{r}=(\:r_0,\:r_1,\:r_2\:)\quad \text{and} \quad \B{\alpha}=(\:\alpha_{\ce{CO2}},\:\alpha_{\ce{H2O}}\:)\ \mathrm{,}
\label{eq:NEfittingparameters}
\end{equation}
Note that even though it is physically necessary to use all of these variables in a fit, the only one of interest in this context is the molecular scaling factor $\alpha_l$ of the molecule $l$ under scrutiny, which in this case is $\alpha_{\ce{CO2}}$), as this alone contains the relevant information about its atmospheric concentration.

This together with \eqref{eq:NEradianceMODEL} clearly fits the criteria for a separable problem, and a conventional VP algorithm from the PORT Mathematical Subroutine Library \citep{Dennis.1981, Fox78} has already been tested by \citet{GimenoGarcia.2011} and validated by \citet{Hochstaffl.2018}. 

\begin{figure}[htbp]
    \centering
    \includegraphics[width=0.95\linewidth]{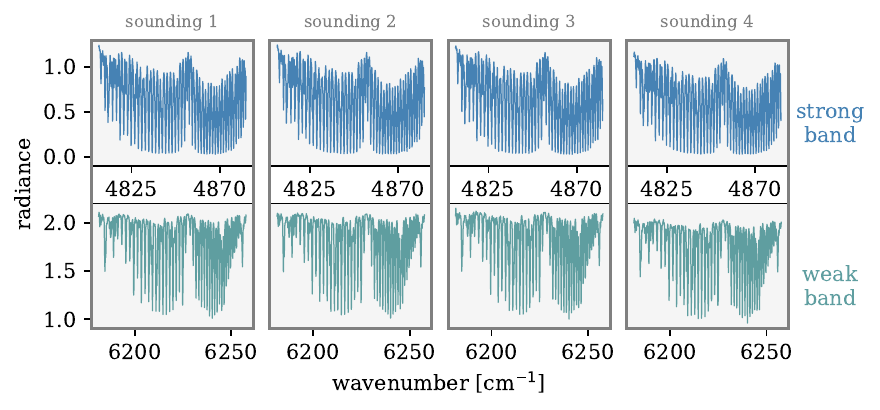}
    \caption{Four exemplary soundings of frame \num{1728} from the OCO-2 level 1b (L1B) measurement product \citep{Crisp}, each displaying radiance spectra in units [$\si{erg /(\s\; \cm\squared\; \steradian\;\cm\tothe{-1}})$] from both the strong and weak bands with \num{809} and \num{651} spectral pixels each.}
    \label{fig:spectra}
\end{figure}

How is this an example of problems with multiple right-handed sides?
Many satellites have sensors that measure radiance simultaneously in several spectral windows, e.g.,\ OCO-2 observes the strong (around \SI{6250}{\cm\tothe{-1}}) and weak absorption bands (around \SI{5000}{\cm\tothe{-1}}) of \ce{CO2} (cf.\ Figure\ \ref{fig:spectra}). Assuming consistent model input, both spectra should deliver the same values for $\B{\alpha}$, but as surface reflectivity varies strongly for different wavelength regions, they each have a specific reflectivity polynomial and, therefore, $\B{r}$. Thus, for every observation, two spectral measurement windows (of different lengths) should be fitted simultaneously.  

This concept of multiple right-hand sides can also be transferred into a spatial dimension: Some molecules, such as carbon dioxide or methane, are very long-lived, so they are distributed relatively homogeneously in the atmosphere. This means that observations from nearby locations should all yield quite similar concentrations. Thus, they might as well be fit for one $\alpha_{\ce{CO2}}$ at once. Note that the assumption made about atmospheric carbon dioxide might not hold for all other absorbing molecules in the observed spectral region, such as \ce{H2O}, which has rather variable concentrations across the globe. However, as their variations are less than that of surface reflectivity, and no physical insight is sought from the fit of the $\alpha_{\ce{H2O}}$ parameter, it can be seen as a mere auxiliary parameter for completing the model and, therefore, be treated as a “constant” nonlinear fitting variable for a group of neighboring spectra. Still, in this case, the reflectivity coefficients, $r_j$, representative of the surface at the place of measurement, are distinct for every geolocation and, therefore, specific to each measured spectrum. Another possible linear model parameter, which is distinct for each spectrum, would be a constant baseline correction added to the model \eqref{eq:NEradianceMODEL}, as suggested by \citet{GimenoGarcia.2011}.

\begin{figure}[htbp]
    \centering
    \includegraphics[width=0.95\linewidth]{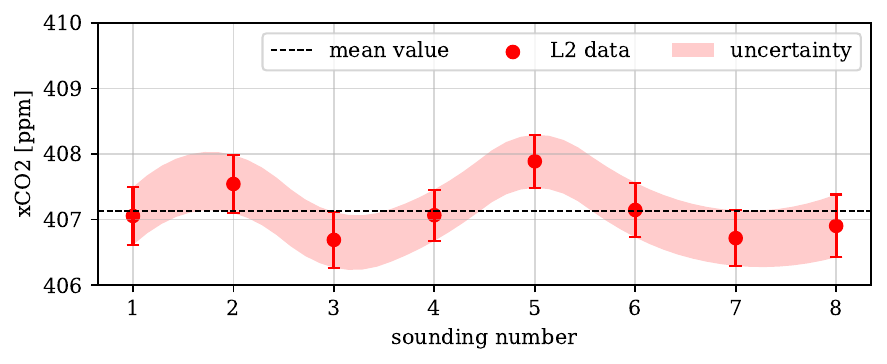}
    \caption{Total column average dry air mole fraction (i.e.,\ “concentration”) of carbon dioxide, denoted as $x\ce{CO2}$ in the unit (parts per million), from the level 2 (L2) retrieval product of \mbox{OCO-2 \citep{L2data, Crisp2}} for one exemplary frame (cf.\ Figure\ \ref{fig:spectra}), including eight soundings.}
    \label{fig:CO2values}
\end{figure}

Finally, it is possible to combine both the spectral and spatial dimensions of multiple RHS fittings in trace gas retrieval. The OCO-2 satellite, for instance, always stores eight observations (“soundings”) in one so-called “frame”, with the spatial coverage not exceeding \SI{24}{\km^2}. The concentration of carbon dioxide can be assumed to fluctuate only minimally within such an area on the globe. Figure\ \ref{fig:CO2values} shows exemplary retrieval results of eight soundings from one OCO-2 frame. Most of the fluctuations in these results (all except for sounding number 5) could be merely due to noise in the measurements, as the mean value stays within the uncertainty for almost all of them. A fit over multiple observations, as proposed, could, therefore, help to constrain the fluctuation to a more reliable retrieval~product.

To summarize, the OCO-2 product \citep{Crisp} allows for multiple RHS fits of a combined \num{8} (spatial) $\cdot$ \num{2} (spectral) $=$ \num{16} datasets (cf.\ Section  \ref{subsec:testsreal}).  
%For the spectral dimension, both bands (strong and weak) have a different number of spectral pixels, meaning different $\B{y}$-sizes, and in the spatial dimension the model requires independent geometry variables and atmospheric state assumptions for every geolocation on the globe, leading to varying model setups $\B{\Phi}_k(\B{\alpha})$, as well. Therefore, this is a suitable example for the modified case of fitting multiple right-hand sides introduced in this paper. 

The following tests were conducted with a \textit{Python} version of BIRRA (Beer infrared retrieval algorithm) \citep{GimenoGarcia.2011}, which has been validated for the SWIR trace gas retrieval of \ce{CO} by \citet{Hochstaffl.2018}. In this code, the Jacobian matrix of the model function \eqref{eq:NEradianceMODEL} is set up analytically for the least squares fit, reducing numerical instabilities. BIRRA is an extension of the radiative transfer model Py4CAtS (Python for Computational Atmospheric Spectroscopy) \citep{Schreier.2019}, which is used to calculate the a priori total optical depths \eqref{eq:NEtotalopticaldepth1} needed in \eqref{eq:NEradianceMODEL}. 

It must be noted that all retrievals conducted with the model described above are only supposed to evaluate the methodology and algorithms and in no way claim to represent full-fledged physical \ce{CO2} retrieval products, such as those by \citet{Crisp2}. Moreover, this technique of fitting multiple spectra measured within a certain spatial distance from each other is only reasonable as long as the assumption holds that, within an order of magnitude of the spatial resolution of the sensor, there are only small to no physically caused fluctuations/gradients in the sought trace gas concentration(s). This is, of course, not the case for localized emitters, such as power plants or biomass-burning events.

\subsection{Tests with Synthetic Data}

The goal of this subsection is to show the conceptual and effective differences of algorithms solving multiple RHS compared to the classical case of solving one. This analysis was conducted on the basis of synthetic spectra. Those are simulated radiance measurements generated with the radiative transfer model Py4CAtS \citep{Schreier.2019}. The benefit of using this in tests is that, in the retrieval (i.e.,\ the fitting process), there is no model error and the exact solution is known. The only deviation from a “perfect” fit is, therefore, controlled by adding noise to the modeled measurements. 

In order to be representative of the later tests with real measurements, the same numbers, sizes, and types (distinct in spatial or spectral dimensions) of datasets were generated as the ones used by OCO-2. Moreover, for consistency, all test retrievals were conducted using the same number of fitting parameters, with $n=\num{3}$ linear ones per dataset and a total of $p=\num{2}$ nonlinear ones (cf.\ vectors in \eqref{eq:NEfittingparameters}). This allowed for the test cases summarized \mbox{in Table \ref{tab:testcases}.}

\begin{table}[htbp]
\caption{Test cases for the analysis of fitting \num{1} (denoted as \textit{single}) or \num{16} (denoted as MRHS) synthetically generated datasets, simultaneously, with both VP and NLS algorithms.}
	\newcolumntype{C}{>{\centering\arraybackslash}X}
	\begin{tabularx}{\textwidth}{m{2cm}<{\centering}m{2cm}<{\centering}m{2cm}<{\centering}m{2.8cm}<{\centering}m{2.8cm}<{\centering}}
	\toprule
\textbf{Name}            & \textbf{\#RHS}  & \textbf{\#fits} & \textbf{\#Parameters per Fit} & \textbf{Tested Algorithm Types}\\ \midrule
\textit{single} & 1  & 16     & $(1\cdot3)+2 = 5$    & VP \& NLS   \\
MRHS   & 16  & 1     & $(16\cdot3)+2 = 50$        & VP \\ \bottomrule
\end{tabularx}
\label{tab:testcases}
\end{table}

For the MRHS case, the Golub--LeVeque VP algorithm (introduced in Section  \ref{subsec:Golub_leVeque Approach}) was used as a representative solver for multiple RHS problems. Tests with synthetic spectra indicated that all mentioned MRHS methods (see Table\ \ref{tab:testcases2}) yielded equal accuracy, confirming the theoretical proof offered by \citet{Golub.1973} that the solutions found by a variable projection solver should be equivalent to those of conventional nonlinear solving methods.

\begin{table}[htbp]
\caption{Test cases for the analysis of real radiance spectra with both VP and NLS algorithms for multiple~RHS.}
	\newcolumntype{C}{>{\centering\arraybackslash}X}
	\begin{tabularx}{\textwidth}{m{3cm}<{\centering}m{4.1cm}<{\centering}m{6cm}<{\centering}}
	\toprule
\textbf{Name} & \textbf{Algorithm}   & \textbf{Reference} \\ \midrule
VP GL & Golub--LeVeque approach & Section \ref{ssec:GOLUB},  \citet{G.H.Golub.1979}  \\ 
VP KM & Kaufman approach & Section \ref{ssec:KAUFMAN},  \citet{Kaufman.1975}  \\
VP naive & Naive approach  & Section \ref{ssec:NAIVE}, \citet{G.H.Golub.1979} \\ 
NLS TRF & Trust Region Reflective & \citet{Branch.1999} \\ 
NLS LM & Levenberg--Marquardt & \citet{More.1978} \\ \bottomrule
\end{tabularx}

\label{tab:testcases2}
\end{table}

For the \textit{single} cases, a classical VP (based on \citet{OLeary.2013}) and a conventional NLS single-RHS solver \citep{SciPyv1.8.0Manual.06.02.2022} (based on \citet{Branch.1999}) were tested. 

First, the fitting precision of the VP MRHS solver was compared to that of \textit{single} solvers using spectra with signal-to-noise  (\textit{SNR}) ratios in the range of \num{20} to \num{500}. One measure of the goodness of fitting results is the relative error,
\begin{equation}
\epsilon = \frac{\alpha_{\text{true}}-\alpha_{\text{fit}}}{\alpha_{\text{true}}},
\end{equation}
compared to the true parameter values. 
%The climatological profile used for the model consists of an a priori number density of \SI{400}{ppm} for \ce{CO2}, which is fairly close to current values, therefore the fitted $\alpha_{\ce{CO2}}$ parameters should also stay close to one.

\begin{figure}[htbp]
    \centering
    \includegraphics[width=0.95\linewidth]{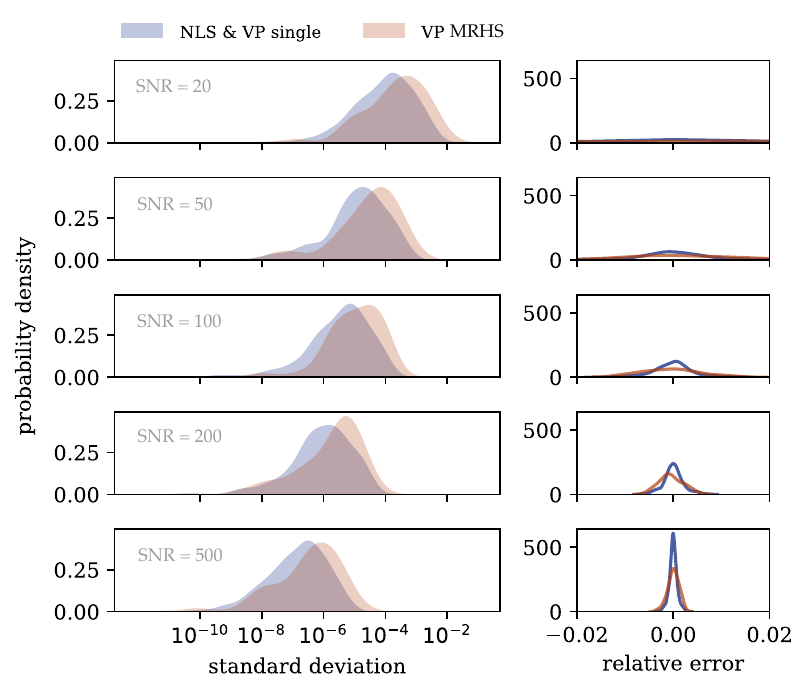}
    \caption{Distribution plots of the relative errors (difference between exact and fitted $\alpha_{\ce{CO2}}$ results) on the right and the corresponding standard deviations from the exact solution on the left for both fitting setups \textit{single} and MRHS with increasing signal-to-noise ratios (\textit{SNRs}).}
    \label{fig:standarddeviations}
\end{figure}

Figure \ref{fig:standarddeviations} shows the distribution of these errors and the corresponding standard deviations for different \textit{SNR} values. The signal-to-noise ratios achieved by satellites, including OCO-2, ranged between approximately \num{200} and \num{800} for the frames used \citep{Crisp}. This broad variation of OCO-2's \textit{SNR} comes from changes in the solar position and varying surface reflectivities across the orbit. As expected, both solvers achieved improved precision for increasing \textit{SNR} values, since a fit becomes more accurate for less noisy data. While both \textit{single} methods (NLS and VP) showed equal performances, the VP MRHS yielded standard deviations that were slightly worse. This trend is also reflected in the distributions of the relative errors, which are always more sharply distributed around zero for the \textit{single} solvers than for the MRHS solver. This behavior is not surprising since a less-dimensional residual vector (coming from the shorter data vector) and fewer unknown parameters most generally leave less freedom in the fit and, therefore, lead to more precise fitting results. To phrase it differently: one can expect that, as the size of the least squares problem increases, the number of possible local minima that the fit can reach will behave accordingly. 

Considering the similarly shaped distributions in Figure\ \ref{fig:standarddeviations} and the fact that VP MRHS seems to improve at the same rate as NLS and VP \textit{single}, MRHS fits can be viewed as equally effective. In particular, at higher \textit{SNR}s, both methodologies achieve deviations from the exact results in such low orders of magnitude that the precisions of their fitted $\alpha_{\ce{CO2}}$ values are very comparable.

Another important measure for the accuracy of a fit is the standard deviation of the residuals $r_i=y_i-\hat{y}_i$ (prediction errors with $\B{y}$: vector of observations, and $\B{\hat{y}}$: fitted model), also known as the sigma of regression, defined as
\begin{equation}
    \sigma = \frac{\norm{\B{r}}}{\sqrt{m-n-p}},
\label{eq:NEsigmaofRegression}
\end{equation}
where the numerator is the norm of the residual vector of the fit and the denominator represents the number of degrees of freedom (the number of data points minus the number of unknown variables). Of course, for the $s$ right-hand side case, the number of linear parameters $n$ becomes $sn$, and the number of data points becomes $ms$ (or $m_1+\ldots+m_s$), respectively.

\begin{figure}[htbp]
    \centering
    \includegraphics[width=0.95\linewidth]{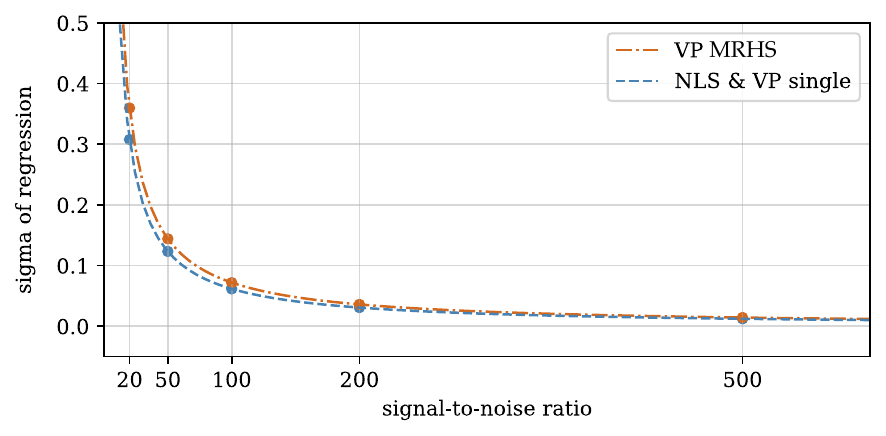}
    \caption{Mean sigma of regression for both \textit{single} and MRHS fits for different noisy spectra. The dashed lines correspond to fitted hyperbolas.}
    \label{fig:regsigma}
\end{figure}

Figure \ref{fig:regsigma} shows the mean sigma of regression of all fits over the \textit{SNR}. Here, the development of \textit{single} and MRHS is similar to that of the errors discussed above. The mean sigmas produced by the MRHS fits are slightly higher than the ones achieved by the \textit{single} solvers. Again, this is intuitive: A \textit{single} solver is able to produce distinct nonlinear parameters for the noisy spectra and, therefore, has more degrees of freedom to mimic the noisy spectra. An MRHS solver, on the other hand, only has one set of nonlinear parameters for all the spectra, leading to an overall larger deviation between the ``observed'' and modeled data. This does not necessarily mean the latter is less accurate. On the contrary, since it is less prone to including the specific noise of the spectra into the fit, MRHS solvers could have a smoothing effect on otherwise fluctuating retrieval results (see Figure\ \ref{fig:CO2values}).

The effect of possible “overfitting” might, however, only be relevant for high noise levels (low \textit{SNR}). As the sigma of regression is proportional to the noisy radiance, which is proportional to $(\B{\num{1}}+\gamma/\text{SNR})$ (with $\B{\gamma}$: normally distributed random value), it decreases by the inverse of the \textit{SNR} value for both the \textit{single} and MRHS fits (see fitted hyperbolas \mbox{in Figure\ \ref{fig:regsigma}).} In this way, as the noise increases, the sigmas of regression become more similar, such that for an \textit{SNR} of \num{200} and higher (representative of OCO-2), the performance differences between the methodologies (\textit{single} and MRHS) disappear. Thus, for reasonably good data, there is no sacrifice in the precision or accuracy of the produced results when fitting multiple RHS simultaneously instead of fitting one by one.

Now that the differences between MRHS solvers and classical \textit{single} solvers are established, we need to analyze which MRHS algorithm (see Section  \ref{subsec:Golub_leVeque Approach}) is the best. 

\subsection{Tests with Real Measurements}
\label{subsec:testsreal}
The goal of this subsection is to assess the performance of the new enhanced VP algorithms for multiple RHS described in Sections \ref{subsec:Extentions to GL} and \ref{subsec:ModificationsforMRHS} (VP naive, VP Golub--LeVeque, VP Kaufman) by comparing them to conventional NLS solvers. 

The \textit{SciPy} function used for the NLS reference approach allows the user to choose between three different nonlinear least squares algorithms (explored thoroughly for a single RHS VP algorithm by \citet{Barligea2022}). In order to better judge the solvers' performances, two of them were used in the tests: the trust region reflective method ('\textit{TRF}') \citep{Branch.1999} and the Levenberg--Marquardt method ('\textit{LM}') \citep{More.1978}. While the former is the most efficient, the latter can be considered the most robust \citep{Barligea2022}, which could be helpful for an increasing number of variables.

In this subsection, an analysis of the test cases listed in Table \ref{tab:testcases2} is conducted.
For the assessment, a set of \num{18} OCO-2 frames (all measured on the 25 of May 2020 on orbit 31366a in the nadir (downward view) acquisition mode just above Australia, with a spacecraft altitude of approximately\ \SI{711}{\km} \citep{L1data}) was used; each included 8 observations in both spectral bands (\mbox{cf.\ Figure\ \ref{fig:spectra})}, within an area of \SI{24}{\km^2} measured along a ground track no wider than \SI{80}{\km}, labeled as cloud-free (no scattering), above land (better reflectivity), and good quality, according to criteria defined by \citet{Crisp}. With those, a few hundred test fits were performed, with the VP and NLS methods (see Table \ref{tab:testcases2}, with varying numbers of RHS ranging from \num{2} to the maximum available number of \num{16} (only even numbers due to the combination of the two spectral bands in one observation). Again, the fits used $n=\num{3}$ linear parameters per spectrum and $p=\num{2}$ nonlinear parameters.

For the evaluation of accuracy, the sigma of regression, the R-Score measure, the confidence bounds of the results, and the fitted residuals, were analyzed. The sigma of regression $\sigma$ defined in Equation \eqref{eq:NEsigmaofRegression} turned out to be equal for all of the tested methods (see the residual analysis below). 

A second statistical quantity is the so-called R-score, defined as 
\begin{equation}
    R =\ \frac{\sum_{i=1}^M (\hat{y}_i - \overline{\B{y}})^2}{\sum_{i=1}^M (y_i - \overline{\B{y}})^2},
\label{eq:NErsquare}
\end{equation}
indicating the amount of variance (the mean of the measurements $\overline{\B{y}}$) accounted for by the fitted model $\hat{\B{y}}$. $M$ means the cumulative size $m_1+\ldots+m_s$ of all the datasets. $R$ must be within $[0,1]$, and the best possible score a fit can achieve would be \num{1}. In the experiments, all of the discussed methods obtained R-scores of approximately\ \num{0.99}. The only difference could be observed for VP GL and VP KM, which had average higher scores by \SI{0.02}{\percent} compared to all other methods, which is negligible. 

In order to calculate the confidence bounds of the retrieval results, the covariance~matrix 
\begin{equation}
    \B{C} = \sigma^2\cdot (\B{H}^{\top}\B{H})^{-1}\ \in\ \R^{(p+ns)\x (p+ns)}
\label{eq:NEcovarianceMatrix}
\end{equation}
needs to be calculated, with $\B{H}$ containing the partial derivatives of the model function, with respect to the $p$ nonlinear and $sn$ linear parameters. For a VP method with multiple right-hand sides, it can be composed as follows:
\begin{equation}
    \B{H} = \begin{pmatrix}
    \frac{\D\B{z}(\B{\alpha})}{\D\B{\alpha}}\big|_{\B{\alpha}=\B{\hat{\alpha}}} & \B{G}(\B{\hat{\alpha}}) 
    \end{pmatrix}\ \in\ \R^{M\x (p+ns)}\mathrm{.}
\label{eq:NEhmatrix}
\end{equation}

Here, the first matrix is the $M\x p$ Jacobian of the purely nonlinear function $\B{z}(\B{\alpha})$ defined in \eqref{eq:THdz(a)/dal} and \eqref{eq:THz(a)NEW} with respect to the nonlinear parameters $\B{\alpha}$, and the second matrix $\B{G}(\B{\hat{\alpha}})\in\R^{M\x ns}$, defined in \eqref{eq:THzeromatrix}, is the Jacobian of all the linear parameters $\B{\beta}$ (\mbox{cf.\ Equation \eqref{eq:THnaiveproblem}). }For a confidence level of \SI{95}{\percent}, one can then calculate the confidence bound(s) (CB) of the retrieved parameters $\hat{\B{x}}$ by
\begin{equation}
    \hat{\B{x}}_{\text{CB}} = q \cdot \sqrt{\text{diag}(\B{C})},
\label{eq:NEconfidencebounds}
\end{equation}
for which $q$ represents the standard normal distribution quantile of $\frac{1}{2}(1-0.95)$, and the diagonal elements of $\B{C}$ are the variances of the estimated parameters $\hat{\B{x}}$. 

\begin{figure}[htbp]
    \centering
    \includegraphics[width=0.95\linewidth]{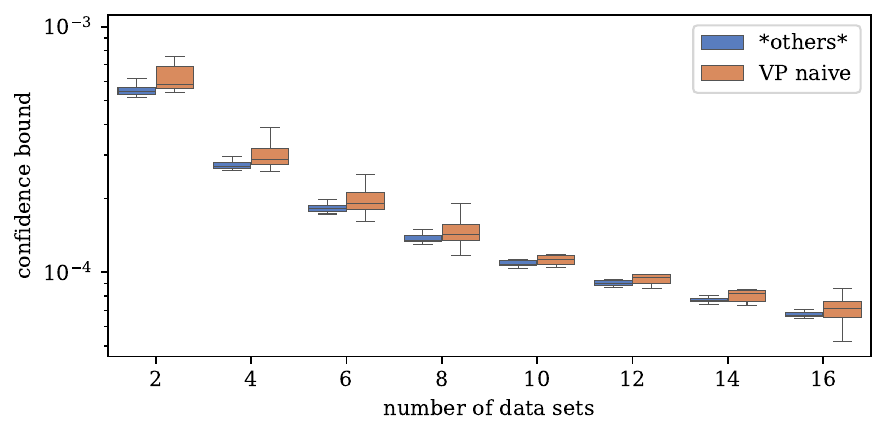}
    \caption{Vertical box plots (with horizontal offset for better distinction) of the confidence bounds of the $\alpha_{\ce{CO2}}$ parameter achieved by each method for different numbers of datasets.}
    \label{fig:bounds}
\end{figure}

Figure \ref{fig:bounds} shows the distribution and mean values of the calculated confidence bounds for the $\alpha_{\ce{CO2}}$ parameter for an increasing number of RHS. While the confidence bounds are already relatively small, they are decreasing for an increasing number of datasets, similar to before with the increasing \textit{SNR} values (see Figure\ \ref{fig:regsigma}). This indicates that more data cause more accurate fitting results. However, a small difference can be observed in Figure\ \ref{fig:bounds} between the naive VP method and the ``others'' (including VP GL, VP KM, NLS TRF, and NLS LM, which all produced the same results). Apparently, the confidence bounds of the results from the naive method, though decreasing, are slightly worse than the rest. This is probably due to the different and more lavishly calculated Jacobian matrix of the naive problem \eqref{eq:THnaiveproblem}. One can, therefore, argue that this is mainly a numerical issue and does not correspond to a lack of accuracy of the VP naive solver. In light of the measures considered above, the tested MRHS solvers all achieved equally accurate fits.

\begin{figure}[htbp]
    \centering
    \includegraphics[width=0.95\linewidth]{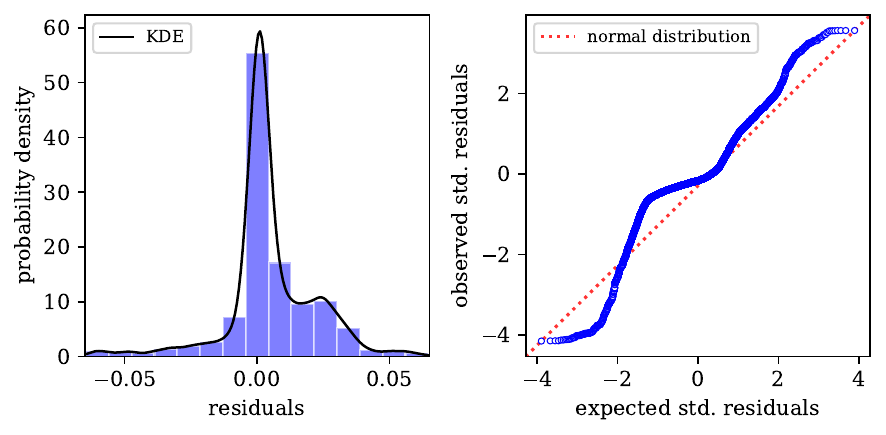}
    \caption{Statistical diagnostics for one exemplary VP fit. The plot on the left shows the distribution of the residuals via a coarse histogram and the continuous kernel density estimate (KDE). The plot on the right is a normal probability plot, displaying the deviation of the residuals from a normal~distribution.}
    \label{fig:residuals}
\end{figure}

This was also confirmed when the residuals of the fits were analyzed, which turned out equally for all methods (cf.\ Table\ \ref{tab:testcases2}). The statistical diagnostics for the residuals of one exemplary VP GL fit (representative of all methods, including NLS) are shown \mbox{in Figure\ \ref{fig:residuals}. }
Ideally, the errors between the fitted model and measurements should be normally distributed. Due to noise and outliers in the spectral data, this distribution may, however, deviate slightly from a normal one. Yet, the fact that the residuals have their highest density around zero indicates that all the algorithms conducted reasonably good~fits.  

As for the robustness, all algorithms yielded convergence rates of \SI{100}{\percent} for decent initial guesses. For a discussion on the impact of bad initial guesses, \mbox{see \citet{OLeary.2013}}, who showed that the VP method ultimately converges more reliably than conventional NLS algorithms. This is mostly due to the fact that the former deal with a reduced nonlinear least squares problem needing only $p$ instead of $p+n$ initial guesses, making the solver a lot more stable.

\begin{figure}[htbp]
    \centering
    \includegraphics[width=0.95\linewidth]{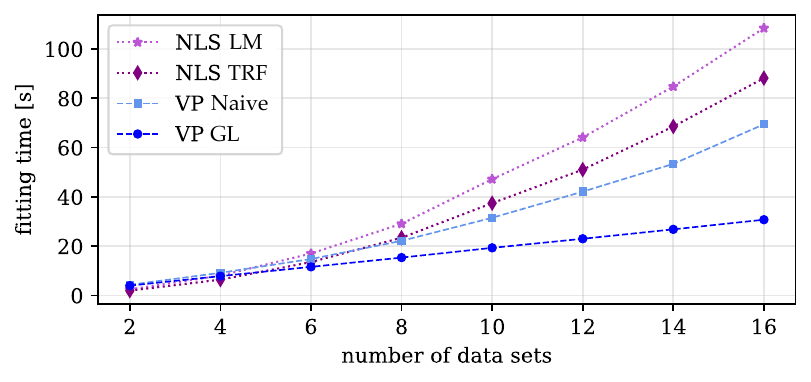}
    \caption{Comparison of the evolution of mean running times of NLS and VP methods for a fit over the number of used datasets.}
    \label{fig:times}
\end{figure}

In the next step, the fitting times of all mentioned methods were analyzed to compare their computational efficiency.
Figure\ \ref{fig:times} shows the mean running times for a fit for an increasing number of datasets. Here, the VP KM method is not shown since its performance is similar to VP GL. 
For fairly small numbers of datasets, the NLS algorithms were faster than the VP methods. This stems from the fact that these algorithms are part of the \textit{SciPy} \mbox{package \citep{SciPy2020}}, which is operationally optimized, whereas the proposed VP code was originally made as a proof of concept and is not yet optimized in the same manner. Still, this scheme changed drastically when more RHS were used. Table \ref{tab:times} shows the exact values for \num{2}, \num{4}, and \num{6} datasets. 

\begin{table}[htbp]
\caption{Fitting times for VP methods compared to NLS methods for small numbers of datasets. The bold numbers mark the smallest mean values within a row.}
\newcolumntype{C}{>{\centering\arraybackslash}X}
	\begin{tabularx}{\textwidth}{m{2cm}<{\centering}m{2cm}<{\centering}m{2cm}<{\centering}m{2.8cm}<{\centering}m{2.8cm}<{\centering}}
	\toprule
        \textbf{\#RHS} & \textbf{VP GL [s]}& \textbf{VP Naive [s]} & \textbf{NLS TRF [s]} & \textbf{NLS LM [s]} \\ \midrule
        \textbf{2} & 3.97 & 4.3 & $\textbf{1.95}$ & 2.46 \\ \midrule
        \textbf{4} & 7.83 & 9.18 & $\textbf{6.40}$ & 8.02 \\ \midrule
        \textbf{6} & $\textbf{11.59}$ & 14.78 & 13.51 & 17.0 \\ \bottomrule
    \end{tabularx}

\label{tab:times}
\end{table}
For six RHS and more, the suggested VP GL algorithm not only becomes significantly faster than the rest, but it is also the only method with fitting times that increase linearly with the number of right-hand sides (see Figure\ \ref{fig:times}), while all the other tested methods exhibit an almost quadratic evolution. This confirms that VP GL and VP KM are the most efficient methods when it comes to dealing with the rising complexity of multiple RHS problems. It also reveals the inferiority of the naive VP method compared to the 'good' VP methods in every test. Even though the naive approach separates the problem and should, therefore, be just as stable as the good approaches, the time it needs for solving also rises quadratically with the number of fitting windows, similar to the slower (and inferior) NLS solvers. This must be due to the increasing size of the block diagonal matrix $\B{G}(\B{\alpha})$ and the resulting extra costs for calculating $\B{G}^{\dagger}(\B{\alpha})$ and the Jacobian $\del\B{P}^{\perp}_{\B{G}}/\del\B{\alpha}$. 

Comparisons of the two ``good'' VP algorithms (VP GL and VP KM) showed that, in all of the above categories, such as robustness or accuracy, the Kaufman approach did equally as well as the Golub--LeVeque one. The only difference could be found in the fitting times, for which the method by \citet{Kaufman.1975}, as predicted, was consistently faster. However, the relative improvements in the running times remained below \SI{1}{\percent} and are, therefore, almost negligible. This confirms the point made by \citet{OLeary.2013} that Kaufman's simplification does not necessarily pose a computational benefit to modern computers anymore. 

\section{Discussion and Conclusions}
\label{sec:Conclusions}
Motivated by a real-world application in atmospheric remote sensing, a variable projection algorithm was extended to multiple right-hand sides of different sizes and nonlinear model setups. A modern \textit{MATLAB} implementation by \citet{OLeary.2013} was translated into \textit{Python} and modified according to the theory presented in this paper. It incorporates the ideas of \citet{G.H.Golub.1979} and \citet{Kaufman.1975} for solving separable nonlinear least squares problems with multiple RHS.

Numerical tests using synthetic data demonstrate that simultaneous fittings over multiple measurements maintain accuracy and precision compared to single dataset solvers, with potential benefits in reducing ``overfitting'' (with noise and outliers affecting the retrieval results) and fluctuations in the results.

A comprehensive comparison with conventional nonlinear least squares solvers using real measurements from NASA's OCO-2 satellite \citep{OCO2} indicated similar accuracy among all algorithms. The most significant finding was that the variable projection methods based on \citet{G.H.Golub.1979} and \citet{Kaufman.1975} significantly outperformed all other methods in computing time, particularly as the number of datasets increased. Thus, these algorithms are deemed more efficient than conventional solvers. Furthermore, our experiments indicate that a popular simplification proposed by \citet{Kaufman.1975} did not yield significant performance improvements.

The algorithm presented in this article proved to be highly effective and efficient. This indicates that the recommended modifications to the original algorithm by \citet{G.H.Golub.1979} preserve its computational advantages. Note, that the benefits arising from a fast solver must always be considered in relation to the computational costs associated with the remaining part of the overall task. In trace gas retrieval applications, the computation time required for the forward model, i.e., radiative transfer \mbox{Equation \eqref{eq:NEradianceMODEL},} can significantly exceed that of solving the inherent least squares problem. Our algorithm, thus, offers the most significant advantage in tasks where the overall performance heavily relies on the fitting process. Consequently, we endorse using this implementation not only for remote sensing, but also for other scientific problems implemented in Python with similar characteristics. The simultaneous fitting of more data can reduce fluctuations in the results, which is highly desirable in some applications.

%%%%%%%%%%%%%%%%%%%%%%%%%%%%%%%%%%%%%%%%%%
\section*{Appendix}
\label{app:A}
In this appendix, derivations for Equations \eqref{eq:INjacobianofVPfunctional} and \eqref{eq:THderivativeQ2} are presented (cf. proof by \citet{Golub.1973}). Beginning with the full Jacobian \eqref{eq:INjacobianofVPfunctional}:
The generalized inverse $\B{\Phi}^{\dagger}=(\B{\Phi}^{\top}\B{\Phi})^{-1}\B{\Phi}^{\top}$ satisfies the identities
\begin{align}
\B{\Phi}\B{\Phi}^{\dagger}\B{\Phi} & =\B{\Phi},\label{eq:AP2}\\
(\B{\Phi}\B{\Phi}^{\dagger})^{\top} & =\B{\Phi}\B{\Phi}^{\dagger}.\label{eq:AP3}
\end{align}

Note that the second identity follows from the fact that the inverse of the symmetric matrix $\B{\Phi}^{T}\B{\Phi}$ is also symmetric. From Equation \eqref{eq:AP3}, we infer that the projector
$\B{P}_{\B{\Phi}}=\B{\Phi}\B{\Phi}^{\dagger}$
is symmetric; from Equation \eqref{eq:AP2}, we infer that 
\begin{align}
\B{P}_{\B{\Phi}}\B{\Phi} & =\B{\Phi}\B{\Phi}^{\dagger}\B{\Phi}=\B{\Phi},\label{eq:AP4}\\
\B{P}_{\B{\Phi}}^{2} & =\B{\Phi}\B{\Phi}^{\dagger}\B{\Phi}\B{\Phi}^{\dagger}=\B{\Phi}\B{\Phi}^{\dagger}=\B{P}_{\B{\Phi}}.\label{eq:AP5}
\end{align}

By means of Equation \eqref{eq:AP4}, we find
\begin{equation}
\del_{l}(\B{P}_{\B{\Phi}}\B{\Phi})=(\del_{l}\B{P}_{\B{\Phi}})\B{\Phi}+\B{P}_{\B{\Phi}}(\del_{l}\B{\Phi})=(\del_{l}\B{\Phi}),\label{eq:AP6}
\end{equation}
implying 
\begin{equation}
(\del_{l}\B{P}_{\B{\Phi}})\B{\Phi}=(\del_{l}\B{\Phi})-\B{P}_{\B{\Phi}}(\del_{l}\B{\Phi})=(\B{I}-\B{P}_{\B{\Phi}})(\del_{l}\B{\Phi})=\B{P}_{\B{\Phi}}^{\perp}(\del_{l}\B{\Phi}),\label{eq:AP7}
\end{equation}
with the short-hand notation $\del_{l}=\del/\del\alpha_{l}$.
Then, using Equation \eqref{eq:AP7} and the relation $\B{P}_{\B{\Phi}}=\B{\Phi}\B{\Phi}^{\dagger}$, we compute the quantity $(\del_{l}\B{P}_{\B{\Phi}})\B{P}_{\B{\Phi}}$ as
\begin{equation}
(\del_{l}\B{P}_{\B{\Phi}})\B{P}_{\B{\Phi}}=(\del_{l}\B{P}_{\B{\Phi}})\B{\Phi}\B{\Phi}^{\dagger}=\B{P}_{\B{\Phi}}^{\perp}(\del_{l}\B{\Phi})\B{\Phi}^{\dagger}.\label{eq:AP8}
\end{equation}

On the other hand, because $\B{P}_{\B{\Phi}}$ is symmetric, $\del_{l}\B{P}_{\B{\Phi}}$ is also symmetric. Therefore, we~have
\begin{equation}
[(\del_{l}\B{P}_{\B{\Phi}})\B{P}_{\B{\Phi}}]^{\top}=\B{P}_{\B{\Phi}}(\del_{l}\B{P}_{\B{\Phi}}),\label{eq:AP9}
\end{equation}
further, in view of Equation \eqref{eq:AP8}, 
\begin{equation}
\B{P}_{\B{\Phi}}(\del_{l}\B{P}_{\B{\Phi}})=[\B{P}_{\B{\Phi}}^{\perp}(\del_{l}\B{\Phi})\B{\Phi}^{\dagger}]^{\top}.\label{eq:AP10}
\end{equation}

Finally, \eqref{eq:AP5}, \eqref{eq:AP8}, and \eqref{eq:AP10} yield 
\begin{equation}
(\del_{l}\B{P}_{\B{\Phi}})=\del_{l}(\B{P}_{\B{\Phi}}^{2})=(\del_{l}\B{P}_{\B{\Phi}})\B{P}_{\B{\Phi}}+\B{P}_{\B{\Phi}}(\del_{l}\B{P}_{\B{\Phi}})=\B{P}_{\B{\Phi}}^{\perp}(\del_{l}\B{\Phi})\B{\Phi}^{\dagger}+[\B{P}_{\B{\Phi}}^{\perp}(\del_{l}\B{\Phi})\B{\Phi}^{\dagger}]^{\top},\label{eq:AP11}
\end{equation}
and the relation $\del_{l}\B{P}_{\B{\Phi}}^{\perp}=-\del_{l}\B{P}_{\B{\Phi}}$
can be used to conclude formula \eqref{eq:INjacobianofVPfunctional} presented in \mbox{Section \ref{subsec:VP}.}

Regarding Equation \eqref{eq:THderivativeQ2}: Using the orthogonality relation $\B{Q}_{2}^{\top}\B{Q}_{1}=\B{0}_{(m-n)\times n}$, we find 
\begin{equation}
\B{Q}_{2}^{\top}\B{\Phi}=\B{Q}_{2}^{\top}\B{Q}_{1}\B{R}_{1}=\B{0}_{(m-n)\times n},\label{eq:AP12}
\end{equation}
so that by differentiating, we have 
\begin{equation}
\del_{l}(\B{Q}_{2}^{\top})\B{\Phi}=-\B{Q}_{2}^{\top}(\del_{l}\B{\Phi}).\label{eq:AP13}
\end{equation}

Multiplying the above Equation from the right by $\B{\Phi}^{\dagger}$ and using the relation $\B{P}_{\B{\Phi}}=\B{\Phi}\B{\Phi}^{\dagger}$, we obtain
\begin{equation}
\del_{l}(\B{Q}_{2}^{\top})\B{P}_{\B{\Phi}}=-\B{Q}_{2}^{\top}(\del_{l}\B{\Phi})\B{\Phi}^{\dagger}.\label{eq:AP14}
\end{equation}

Now, using the identity
\begin{equation}
\B{Q}_{2}^{\top}\B{P}_{\B{\Phi}}=\B{Q}_{2}^{\top}(\begin{array}{cc}
\B{Q}_{1} & \B{Q}_{2}\end{array})\left(\begin{array}{c}
\B{Q}_{1}^{\mathrm{\top}}\\
\B{Q}_{2}^{\top}
\end{array}\right)=(\begin{array}{cc}
\B{0}_{(m-n)\times n} & \B{I}_{(m-n)\times(m-n)}\end{array})\left(\begin{array}{c}
\B{Q}_{1}^{\mathrm{\top}}\\
\B{Q}_{2}^{\top}
\end{array}\right)=\B{Q}_{2}^{\top},\label{eq:AP15}
\end{equation}
yields
\begin{equation}
\del_{l}(\B{Q}_{2}^{\top}\B{P}_{\B{\Phi}})=\del_{l}(\B{Q}_{2}^{\top})\B{P}_{\B{\Phi}}+\B{Q}_{2}^{\top}(\del_{l}\B{P}_{\B{\Phi}})=\del_{l}(\B{Q}_{2}^{\top}).\label{eq:AP16}
\end{equation}

We infer that (cf. Equation \eqref{eq:AP14})
\begin{equation}
\del_{l}(\B{Q}_{2}^{\top})=-\B{Q}_{2}^{\top}(\del_{l}\B{\Phi})\B{\Phi}^{\dagger}+\B{Q}_{2}^{\top}(\del_{l}\B{P}_{\B{\Phi}}),\label{eq:AP17}
\end{equation}
where $\del_{l}\B{P}_{\B{\Phi}}$ is given by \eqref{eq:AP11}. However, \citet{Kaufman.1975} suggested and mathematically justified neglecting the second term in the above Equation. The resulting approximation
is Equation \eqref{eq:THderivativeQ2}, which works sufficiently well in practice.

%%%%%%%%%%%%%%%%%%%%%%%%%%%%%%%%%%%%%%%%%%
\section*{Data Availability}
The \textit{Python} code for the presented variable projection algorithm, as well as examples, are available at \url{https://atmos.eoc.dlr.de/tools/varPro/index.html} (accessed on 16 May 2023). All data used for the tests with synthetic spectra can be reproduced with the radiative transfer model Py4CAtS (\url{https://atmos.eoc.dlr.de/tools/Py4CAtS/} (accessed on 16 May 2023)) \citep{Schreier.2019}.
Data concerning NASA's OCO-2 satellite are publically available at  \url{https://disc.gsfc.nasa.gov/} (accessed on 16 May 2023), provided by the Goddard Earth Sciences Data and Information Services Center (GES DISC) \citep{L1data, L2data}. 
Any data arising from the numerical experiments that were used in the evaluations will be made available upon request. 

\bibliographystyle{unsrtnat}
\bibliography{reference} 

\end{document}